\documentclass{article}
\newcommand{\R}{I\!\! R}
\newcommand{\N}{I\!\! N}

\begin{document}

\centerline{\bf Affine Anosov diffeomorphims of affine manifolds.}

\bigskip

 by Aristide Tsemo

The Abdus Salam International Centre for Theoretical Physics

Strada Costiera, 11

Trieste Italy.

Email tsemo@ictp.trieste.it

\bigskip

\centerline{\bf Abstract.}

{\it In this paper, we show that a compact affine manifold endowed
with an affine Anosov transformation map is finitely covered by a
complete affine nilmanifold. This answers a conjecture of Franks
for affine manifolds.}

\bigskip
\bigskip

\centerline{0. Introduction.}

\bigskip

An $n-$affine manifold $(M,\nabla_M)$ is an $n-$differentiable
manifold $M$ endowed with  a locally flat connection $\nabla_M$,
that is a connection $\nabla_M$ whose curvature and torsion forms
vanish identically. The connection $\nabla_M$ defines on $M$ an
atlas (affine) whose transition functions are affine maps. The
connection $\nabla_M$ pulls-back on the universal cover $\hat M$
of $M$, and defines on it a locally flat connection ${\nabla_{\hat
M}}$. The affine structure of $\hat M$ is defined by a local
diffeomorphism $D_M:\hat M\rightarrow {\R}^n$ called the
developing map. The developing map gives rise to a representation
$h_M:\pi_1(M)\rightarrow Aff({\R}^n)$ called the holonomy. The
linear part $L(h_M)$ of $h_M$ is the linear holonomy. The affine
manifold $(M,\nabla_M)$ is complete if and only if $D_M$ is a
diffeomorphism. This means also that the connection $\nabla_M$ is
complete.

A diffeomorphism $f$ of $M$ is an Anosov diffeomorphism if and
only if there is a norm $\mid\mid$ $\mid\mid$ associated to a
riemannian metric $<,>$, a real number $0<\lambda<1$ such that the
tangent bundle $TM$ of $M$ is the direct summand of two bundles
$TM^s$ and $TM^u$ called respectively the stable and the unstable
bundle such that:
$$
\mid\mid df^m(v)\mid\mid \leq c\lambda^m\mid\mid v\mid\mid, v\in
TM^s
$$

$$
\mid\mid df^m(w)\mid\mid \geq c\lambda^{-m}\mid\mid w\mid\mid,
w\in TM^u;
$$
where $d$ is the usual differential, $c$ a positive real number,
and $m$ a positive integer.

The stable distribution $TM^s$ (resp. the unstable distribution
$TM^u$) is tangent to a topological foliation ${\cal F}^s$ (resp.
${\cal F}^u$).

The property for a diffeomorphism to be  Anosov is independent of
the choice of the riemannian metric if $M$ is compact. In this
case we can suppose that $c$ is $1$.

\medskip

A Franks's conjecture asserts that an Anosov diffeomorphism $f$ of
a compact manifold $M$ is $C^0-$conjugated to an hyperbolic
infranilautomorphism. This conjecture is proved in [1] with the
assumptions that $f$ is topologically transitive it's stable and
unstable foliations are $C^{\infty}$, it preserves a symplectic
form or a connection.

The goal of this paper is to characterize compact affine manifolds
endowed with affine Anosov transformations. More precisely, we
show:

\bigskip

{\bf Theorem 1.} {\it Let $(M,\nabla_M)$ be a compact affine
manifold, and $f$ an affine Anosov transformation of $M$, then
$(M,\nabla_M)$ is finitely covered by a complete affine
nilmanifold.}

\bigskip

{\bf 1. The Proof of the main theorem.}

\bigskip

The main goal of this part is to show theorem $1$. In the sequel,
$(M,\nabla_M)$ will be an $n-$compact affine manifold endowed with
an affine Anosov diffeomorphism $f$. The stable foliation ${\cal
F}^s$ (resp. the unstable foliation ${\cal F}^u$) pulls-back on
the universal cover $\hat M$ to a foliation $\hat {\cal F}^s$,
(resp. $\hat {\cal F}^u$).

Let $(U,\phi)$ be an affine chart of $M$, the restriction of a
riemannian metric of $M$ to $U$ will be said an euclidean metric
adapted to the affine structure of $U$ if its restriction to $U$
is the pulls-back by $\phi$ on $U$, of the restriction to
$\phi(U)$ of an euclidean metric of ${\R}^n$.

\bigskip

{\bf Proposition 1.1.} {\it The stable and the unstable
distributions of $f$ define on $(M,\nabla_M)$ foliations whose
leaves are immersed affine submanifolds.}

\medskip

{\bf Proof.} Let $x$ be an element of $M$, $\mid\mid$ $\mid\mid$ a
norm associated to a riemannian metric $<,>$ of $M$. Let $v$ be an
element of $T_xM^s$ the subspace of $T_xM$ tangent to ${\cal
F}^s$, we have $\mid\mid df^m(v)\mid\mid\leq \lambda^m\mid\mid
v\mid\mid$ where $0<\lambda<1$, and $m\in {\N}$. Let $(U,\phi)$ be
an affine chart which contains a point of accumulation of the
sequence $(f^p(x))_{p\in{\N}}$. We can suppose that the
restriction of $<,>$ to $U$ is euclidean and adapted to the affine
structure. Let $p>p'$ such that $f^p(x)$ and $f^{p'}(x)$ are
elements of $U$, and $v\in T_{f^{p'}(x)}M^s$, we have $\mid\mid
df^{p-p'}(v)\mid\mid\leq \lambda^{p-p'}\mid\mid v\mid\mid$. This
implies that for every element $y=f^{p'}(x)+w, w\in
T_{f^{p'}(x)}M^s$ and $f^{p-p'}(y)\in U$ is an element of the
stable leaf of $f^{p'}(x)$ since the distance between $f^{n_q}(x)$
and $f^{n_q}(y)$ goes to zero for a subsequence $n_q>p$ such that
$f^{n_q}(x)$ is an element of $U$. We deduce that this leaf is an
immersed submanifold. The result for the unstable foliation is
deduced considering $f^{-1}$.

\bigskip

{\bf Proposition 1.2.} {\it The leaves of ${\cal F}^s$ are
geodesically complete for the affine structure.}

\medskip

{\bf Proof.} Let $x$ be an element of $M$, as $M$ is supposed to
be compact, the sequence $(f^m(x))_{m\in {\N}}$ has a point of
accumulation $y$. Let $U_y$ be an open set containing $y$ such
that there is a strictly positive number $r$, such that for every
$z\in U_y$, $v\in T_zM$ whose norm is less than $r$ for a given
riemannian metric, a (affine) geodesic from $x$ whose derivative
at $0$ is $v$ is defined at $1$. Let $w$ be an element of $T_xM^s$
such that $\mid\mid df^p(w)\mid\mid\leq r$.
 Without loss generality, we can
suppose that $f^p(x)\in U_y$, the (affine) geodesic from $f^p(x)$
whose derivative is $df^p(w)$ at $0$ is defined at $1$. This
implies that the geodesic from $x$ which derivative at $0$ is $w$
is defined at $1$, since $f^{-1}$ is an affine map. This shows the
result.

\bigskip

It is a well-known fact that an Anosov's diffeomorphism of a
compact manifold as a periodic fixed point. We will assume up to
change $f$ by an iterated that f has a fixed point $x$. There
exists also a map $F:\hat M\rightarrow \hat M$ over $f$ which
fixed the element $\hat x$ over $x$.

\bigskip

{\bf Proposition 1.3.} {\it Let $\hat y$ and $\hat z$ be two
elements of $\hat{\cal F}^u_{\hat t }$, where $t$ is an element of
$\hat M$, then the images of $\hat{\cal F}^s_{\cal y}$ and
$\hat{\cal F}^s_{\hat z}$ by the developing map are parallel
affine subspaces.}

\medskip

{\bf Proof.} Let $\hat y$ be an element of $\hat{\cal F}^u_{\hat
t}$, it is sufficient to show that $D(\hat{\cal F}^s_{\hat t})$
and $D(\hat{\cal F}^s_{\hat y})$ are parallel affine subspaces.

We know that the tangent bundle of a simply connected affine
manifold is trivial. We have $T\hat M=\hat M\times T_{\hat t}\hat
M$. Let $<,>'$ be a riemannian metric on $M$ which pulls back on
$\hat M$ is $<,>$. The map $F$ is Anosov relatively to $<,>$.
Without loss generality, we can assume that the distributions
tangent to $\hat{\cal F}^s$ and $\hat{\cal F}^u$ are orthogonal.

Let $w$ be a vector of $T\hat M_{\hat y}$ tangent to $\hat{\cal
F}^s_{\hat y}$. We can write $w=s+u$ where $s$ is a vector of
$T\hat M_{\hat t}$ tangent to $\hat{\cal F}^s$ and $u$ is a vector
of $T\hat M_{\hat t}$ tangent to $\hat{\cal F}^u$. The vector $u$
is also an element of $T\hat M_{\hat y}$ tangent to $\hat{\cal
F}^u_{\hat y}$ since the unstable foliation is affine (we can
identify the projection on the second factor of $T\hat M\times
T\hat M_{\hat t}$ of vectors tangent to $\hat t$ and $\hat y$
since $T\hat M$ is trivial). This implies that $\mid\mid
dF^n_{\hat y}(u)\mid\mid\geq \lambda^{-n}\mid\mid u\mid\mid$ for
$0<\lambda<1$. But on the other hand we have $\mid\mid dF^n_{\hat
y}(s+u)\mid\mid\leq \lambda^n\mid\mid (s+u)\mid\mid$, which
implies that the limit of $\mid\mid dF^n_{\hat y}(s+u)\mid\mid$ is
zero when $n$ goes to infinity. We deduce that $u=0$ since we have
supposed that the distributions tangent to $\hat{\cal F}^s$ and
$\hat{\cal F}^u$ are orthogonal, so $\mid\mid dF_{\hat
y}(s+u)\mid\mid=\mid\mid dF_{\hat y}(s)\mid\mid+\mid\mid dF_{\hat
y}(u)\mid\mid$. This implies the result.

\bigskip

 The images by $D_M$ of the leaf
$\hat{\cal F}^s_{\hat t}$ of $\hat {\cal F}^s$ and  $\hat{\cal
F}^u_{\hat t}$ of $\hat {\cal F}^u$  are affine subspaces of
${\R}^n$ whose direction are direct summand of ${\R}^n$. Since for
every element  $z$ of    $\hat{\cal F}^s_{\hat t}$ the leaf of
${\cal F}^u$ passing by  $z$ is complete, we deduce that the
developing map is surjective.

\bigskip

{\bf Proposition 1.4.} {\it The affine manifold $(M,\nabla_M)$ is
complete.}

\medskip

{\bf Proof.} Let $\hat t$ be an element of $\hat M$, and $\hat
E_{\hat t}$ the set of elements  $y$ of $\hat M$ such that there
is and element $z$ in $\hat{\cal F}^u_{\hat t}$ such that $y$ is
an element of $\hat{\cal F}^s_z$. The image of $\hat E_{\hat t}$
by $D_M$ is ${\R}^n$, and the restriction of $D_M$ to $\hat
E_{\hat t}$ is injective. The set $\{\hat E_{\hat t},\hat t \in
\hat M\}$ is a partition of $\hat M$ by disjoint open sets. It has
only one element since $\hat M$ is connected. We deduce that
$(M,\nabla_M)$ is complete.

\bigskip

{\bf Remark.}

\medskip

The existence of an affine Anosov transformation $\hat f$ on the
universal cover of a compact affine manifold $(M,\nabla_M)$ which
pulls forward onto a diffeomorphism $f$ of $M$ does not implies
that $f$ is an Anosov diffeomorphism as shows the following
example:

Let $Hopf(n)$ be the quotient of ${\R}^n/\{0\}$ by an homothetie
$h_\lambda$ whose ratio $\lambda$ is such that $0<\lambda<1$. It
is a compact affine manifold. Every homothetie $h_c$ which ratio
$c$ is positive and  different from $1$ and $\lambda$ is an Anosov
diffeomorphism of ${\R}^n$ endowed with an euclidean metric. But
the pulls forward of $h_c$  is an isometry of $Hopf(n)$ endowed
with the pulls forward of the riemannian metric of ${\R}^n/\{0\}$
defined as follows:
$$
{1\over\mid\mid x\mid\mid^2}<u,v>
$$
where $x$ is an element of ${\R}^n/\{0\}$, and $u$, $v$ are
elements of its tangent space.

\bigskip

{\bf Proof of theorem 1.}

\bigskip

First we show that the module of the eigenvalues of the elements
of the linear holonomy of $(M,\nabla_M)$ is $1$.

Let $(A,a)$ be an element of $\pi_1(M)$ such that $A$ has an
eigenvector $u$ associated to an eigenvalue $b$ (which may be a
complex number) whose norm is different from $1$.

It is a well-known fact that an Anosov diffeomorphism of a compact
manifold has a fixed periodic point. Without loss of generality,
we will assume that up to change $f$ by an iterated that $f$ has a
fixed point $x$. This implies that up to a change of coordinates,
there exists on ${\R}^n$ a linear map $F$ over $f$.

Put $x=p(0)$ where $p$ is the covering map and consider a
riemannian metric $<,>'$ on $M$ whose restriction on an affine
neighborhood $N$ of $x$ is euclidean adapted to the affine
structure. Expressing on $N$ the fact that $F$ is an Anosov
diffeomorphism using the metric $<,>'$, one obtains that
${\R}^n=U\oplus V$, where $U$ and $V$ are two subvectors spaces
such that there exists a number $0<\lambda<1$ such that
$$
\mid\mid F^n(u)\mid\mid\leq \lambda^n\mid\mid u\mid\mid, u\in U,
$$
$$
\mid\mid F^n(v)\mid\mid\geq \lambda^{-n}\mid\mid u\mid\mid, v\in V
$$
where $u$, and $v$ are respectively elements of $U$ and $V$ and
$\mid\mid$ $\mid\mid$ is a norm associated to an euclidean metric
$<,>$ of ${\R}^n$. The subvectors spaces $U$ and $V$ pulls forward
respectively onto $T_xM^s$ and $T_xM^u$. We will assume that they
are orthogonal with respect to $<,>$. The vectors spaces $U$ and
$V$ are stable by the linear holonomy (see proposition 1.3)

Put $u=u_1+u_2$, where $u_1$ and $u_2$ are respectively elements
of $U'$ and $V'$ the complexified vectors spaces respectively
associated to $U$ and $V$.

Without restrict the generality, one can assume that after
eventually having changed $\gamma$ by $\gamma^{-1}$ and (or) $f$
by $f^{-1}$ that $u_1$ is not zero and that the norm of $b$ is
strictly superior to $1$.

Let $q$ be a positive integer, consider the smallest positive
integer $n_q$ such that $F^{n_q}\circ\gamma^q$ has a fixed point
$\hat m_q$. The integer $n_q$ exists since $A$ preserves $U$ and
$V$. Let $y$ be a point of accumulation of the sequence $p(\hat
m_q)$. Up to replace $(\hat m_q)$ by a subsequence and replace
$\hat m_q$ by another element over $m_q$, we can suppose that the
sequence $(\hat m_q)$ converges to $\hat y$ over $y$. The linear
part of the elements over $f^{n_q}$ which fixe the elements over
$m_q$ have the same eigenvalues since they are conjugated.
Consider a riemannian metric of $M$ which pulls back on ${\R}^n$
coincides on a neighborhood of $\hat y$ with an euclidean metric.

$1.$ The sequence $(n_q)$ is bounded.

The restriction of the linear part of the element of $f^{n_q}$
conjugated to $F^{n_q}\circ \gamma^q$ which fixe $\hat m_q$ at $U$
is contractant  (for the last euclidean metric) with ratio
$\lambda^{n_q}$ ($0<\lambda<1$). This is not possible  since the
restriction of the linear part of this element to $U'$ has the
same eigenvalues than the restriction of $F^{n_q}\circ A^q$ to
$U'$, and the limit when $q$ goes to infinity of the norm of
$F^{n_q}\circ A^q(u_1)$ for the hermitian metric associated to the
last euclidean metric is infinity, since the sequence $n_q$ is
bounded and the norm of $b>1$.

\medskip

$2$. The sequence $(n_q)$ is not bounded.

Up to change $(n_q)$ by a subsequence, we can suppose that $(n_q)$
goes to infinity. The linear map $F^{n_q-1}\circ\gamma^q$ does not
have a fixed point. Its linear part has the eigenvalue $1$
associated to the eigenvector $v_q$. Write $v_q=v_{1q}+v_{2q}$
where $v_{1q}$ and $v_{2q}$ are respectively element of $U$ and
$V$. If $v_{1q}$ is not zero, then $\mid\mid F^{n_q}\circ
A^q(v_{1q})\mid\mid=\mid\mid F(v_{1q})\mid\mid$ (the norm
considered is the precedent euclidean norm); the restriction of
$F^{n_q}\circ A^q$ cannot be contractant with ratio
$\lambda^{n_q}$ for $q$ enough big since the sequence $(n_q)$ goes
to infinity.

If $v_{2q}$ is not zero, then $\mid\mid F^{n_q}\circ
A^q(v_{2q})\mid\mid=\mid\mid F(v_{2q})\mid\mid$, the restriction
of $F^{n_q}\circ A^q$ to $V$ cannot be dilatant (for the last
euclidean metric) with ratio $\lambda^{-n_q}$ for $q$ enough big
since $n_q$ goes to infinity.
 There is a contradiction since the eigenvalues of the
restriction of $F^{n_q}\circ A^q$ to $U'$ (resp to $V'$) coincide
with the eigenvalues of the restriction to $U'$ of the linear part
of the conjugated to $F^{n_q}\circ \gamma^q$ which fixes $\hat
m_q$ (resp. the eigenvalues of the restriction to $V'$ of the
linear part of the element conjugated to $F^{n_q}\circ \gamma^q$
which fixes $\hat m_q$). We deduce that the eigenvalues of $A$ are
roots of unity.

We deduce from the first paragraph of  [4] p.6 that
 $\pi_1(M)$ has a
subgroup of finite index $G$ such that the eigenvalues of the
linear part of the elements of $G$ are $1$. We deduce from [3]
that the quotient of ${\R}^n$ by $G$ is a complete affine
nilmanifold.

\bigskip
\bigskip

{\bf Remark.}

Let $(M',\nabla')$ be the finite cover of $(M,\nabla)$ which
fundamental group is $G$, and $f'$ the pulls-back of $f$ to $M'$.
The map $H^n(M,{\R})\rightarrow H^n(M,{\R})$, $\alpha\rightarrow
(f'^2)^*\alpha$ is the identity since $f'^2$ is a diffeomorphism
which preserves the orientation. Let $\omega$ be the parallel form
of $(M',\nabla')$, we deduce that $(f'^2)$ preserves the parallel
volume form of $(M',\nabla')$; then a well-known result of Anosov
implies that $f'^2$ is ergodic and its periodic points are dense.

\bigskip
\bigskip

\centerline{\bf Bibliography.}

\bigskip

[1] Benoist Y., Labourie F., Sur les diff\'eomorphismes d'Anosov
affines \`a feuilles stables et instables diff\'erentiables,
Invent. Math. 111 (1993) 285-308.

[2] Benoist Y., Foulon P., Labourie F., Flots d'Anosov \`a
distributions stable et instable diff\'erentiables, J.  Amer.
Math. Soc. 5 (1992)  33-74.

[3] Fried D., Goldman, W., Hirsch M., Affine manifolds with
nilpotent holonomy, Comment. Math. helv. 56 (1981) 487-523.

[4] Goldman W., Hirsch M., A generalization of Bieberbach's
theorem, Invent. Math. 65 (1981) 1-11.

[5] Koszul J., Vari\'et\'es localement plates et convexit\'e,
Osaka J. Math. 2 (1965) 285-290.

[6] Smale S., Differentiable dynamical systems, Bull. Amer. Math.
Soc. 73 (1967) 747-817.

\end{document}